\documentclass[11pt, reqno]{amsart}
\usepackage{amsmath, amsthm, amscd, amsfonts, amssymb, graphicx, color}
\usepackage[bookmarksnumbered, colorlinks, plainpages]{hyperref}
\usepackage{enumitem}

\makeatletter \oddsidemargin.9375in \evensidemargin \oddsidemargin
\def\authorsaddresses#1{\dedicatory{#1}}
\newtheorem{theorem}{Theorem}[section]
\newtheorem{lemma}[theorem]{Lemma}

\newtheorem{corollary}[theorem]{Corollary}
\theoremstyle{definition}

\newtheorem{example}[theorem]{Example}

\theoremstyle{remark}
\newtheorem{remark}[theorem]{Remark}
\theoremstyle{approach}

\numberwithin{equation}{section}

\begin{document}
\setcounter{page}{1}


\title[ON THE ARENS  REGULARITY OF THE HERZ ALGEBRAS
$A_p (G)$ 
 ]{ON THE ARENS  REGULARITY OF THE HERZ ALGEBRAS
$A_p (G)$ 
}

\author[H. G. AMINI AND A.REJALI]{H. G. AMINI AND A.REJALI}
\subjclass[2000]{ Mathematics Subject Classification. Primary 43A15; Secondary 46AJ10.}

\keywords{Locally compact groups, Herz algebra, Fourier and Fourier-Stieltjes
algebra, amenability.}

\begin{abstract}
Let $G$ be a locally compact group, $A_p (G)$ be the Herz algebra of $G$ associated with  $1 <p< \infty$. We show that $A_p (G)$ is Arens regular if and only if $G$ is a discrete group and for each countable subgroup $H$ of $G$, $A_p (H)$ is Arens regular. In the case $G$ is a countable discrete group we investigate the relations between Arens regularity of $A_p (G)$ and the iterated limit condition. We consider the problem of Arens regularity of $l^1 (G)$ as a subspace of $A_p (G)$. A few related results when the unit ball of $(l^1 (G),.,A_p(G))$ is bounded under $\|.\|_1$-norm are also determined.
\end{abstract}

\maketitle


\section{Introduction}
Let $G$ be a locally compact group. For $1 <p< \infty$, let $A_p (G)$ denote the linear subspace of $C_0 (G)$ consisting of all functions of the form
\begin{equation*}
u=\sum_{i=1}^{\infty}g_i*f^{\vee}_i
\end{equation*}
where $f_i\in L_p(G)$, $g_i\in L_q(G)$, $\frac{1}{p}+\frac{1}{q}=1$ and $\sum_{i=1}^{\infty}\|f_i\|_p\|g_i\|_q<\infty$, $f^{\vee}(x)=f(x^{-1})$ for $x\in G$. $A_p(G)$ is a commutative Banach algebra with respect to the pointwise multiplication and the norm,
\begin{equation*}
\|u\|_{A_p(G)}=\inf\{\sum_{i=1}^{\infty}\|f_i\|_p\|g_i\|_q: u=\sum_{i=1}^{\infty}g_i*f^{\vee}_i\}.
\end{equation*}
When $p=2$, $A_2(G) = A(G)$ is the Fourier algebra of $G$ as introduced by Eymard \cite{9}. The algebras $A_p(G)$ was introduced and studied by Herz \cite{13}. Let $P M_p(G)$ denote the closure of $L_1(G)$, considered as an algebra of convolution operators on $L_p(G)$, with respect to the weak operator topology in the bounded operators on $L_p(G)$, denoted by $B(L_p(G))$. The space $P M_p(G)$ can be identified with the dual of $A_p(G)$ for each $1 <p< \infty$ \cite{6}. \\
\indent
Arens regularity of $A_p(G)$, has been studied by B. Forrest \cite{10}. It is shown, for example, that if the algebra $A_p(G)$ is Arens regular then the group $G$ is discrete. Furthermore, as shown in \cite{15}, Proposition 5.3, for $p = 2$ and $G$ amenable, $A_p(G) = A(G)$ is Arens regular if and only if $G$ is finite. In this paper, we shall study some problems related to Arens regularity of $A_p(G)$.\\
\indent
Notations and definitions which are necessary in the sequel are gathered in \S2. In \S3, it ts shown that the Arens regularity of $A_p(G)$ can be reduced when G is discrete and finitely generated. Furthermore, we show that $A_p(G)$ is Arens regular if and only if $A_p(H)$ is Arens regular, for any countable subgroup $H$ of $G$. In the case $G$ is a discrete group, we consider the iterated limit condition for the Arens regularity of $A_p(G)$. \\
\indent
In \cite{3} J. W. Baker and A. Rejali considered the Arens regularity of the algebras $l_1(S)$ where $S$ is a discrete semigroup. Similarly, in Theorem 3.6 we consider the Arens regularity of $l_1(G)$ under various multiplications and norms. Furthermore, we show that $(l_1(G), ., \|.\|A_p(G))$ is Arens regular if and only if $A_p(G)$ is Arens regular. In Theorem 3.7 we also show that if the unit ball of $(l_1(G), ., \|.\|A_p(G))$ is bounded under $\|.\|_1$, then $G$ has no infinite Abelian subgroup.
\section{Preliminaries and some notations}
Let $G$ be a locally compact group with a fixed left Haar measure $\lambda$, which is also denoted by $dx$. The usual $L_p(G)\,  (1 \leq p \leq\infty)$ are defined with respect to this measure $\lambda$. By $C(G)$ we denote the Banach space of bounded continuous complex valued functions endowed with the supremum norm. By $C_0(G)$ (resp. $C_{00}(G)$) we denote the subspace of $C(G)$ consisting of the functions on $G$ vanishing at infinity (resp. with compact support). \\
\indent
For a Banach space $X$, we denote by $X^*$ and $X^{**}$ its first and second continuous duals. We regard X as naturally embedded into $X^{**}$. For $x\in X$ and $x^*\in X^*$, by $< x, x^* >$ (or by $< x^*, x >$) we denote the natural duality between $X$ and $X^*$. \\
\indent
Let $A$ be a Banach algebra. For $u\in A$ and $T\in A^*$ by $uT$ (resp. $Tu$) we denote the element of $A^*$ defined by $< uT, v >=< T, vu >$ (resp. $< Tu, v >=< T, uv >$), for $v$ in $A$. When $A$ is commutative it is clear that $uT = Tu$. Then $A^{**}$ can be given two multiplications that extend the multiplication of $A$ and for which $A^{**}$ becomes a Banach algebra. When these two multiplications coincide on $A^{**}$, the algebra $A$ is said to be ‘Arens regular’. Details of the constructions can be found in many places, including the book \cite{4} and the papers \cite{1, 2, 7, 17}. \\
\indent
Let $A$ be an arbitrary Banach algebra and $A_1$ be its closed unit ball. We call $T\in A^*$ weakly almost periodic if $O(T ) = \{T u : u \in A_1\}$ is relatively weakly compact in $A^*$. By $wap(A)$ we denote the linear subspace of $A^*$ consisting of all weakly almost periodic functionals on $A$. The space $wap(A)$ is a closed subspace of $A^*$ and by Theorem 1 of \cite{7}, $A$ is Arens regular if and only if $wap(A) = A^*$ if and only if for each sequences $(a_n), (b_m)$ in $A_1$ and $T\in A^*$, the iterated limits
\begin{equation*}
\lim_n \lim_m T(a_nb_m),\qquad	\lim_m \lim_n T (a_nb_m)
\end{equation*}
are equal, whenever both exist.
\section{The Arens regularity and iterated limit conditions}
Forrest \cite{10} showed that $G$ is discrete, whenever $A_p(G)$ is Arens regular. In the following, we obtain a criteria for the Arens regularity of $A_p(G)$.
\begin{theorem}
 Let $G$ be a locally compact group. Then $A_p(G)$ is Arens regular if and only if $G$ is a discrete group and for each countable subgroup $H$ of $G$, $A_p(H)$ is Arens regular.
\end{theorem}
{\bf Proof.}
If $A_p(G)$ is Arens regular, then by \cite[Theorem 3.2]{10} $G$ is discrete and for each subgroup $H$ of $G$, $A_p(H)$ is Arens regular \cite[Lemma 3.1]{10}.\\
\indent
Conversely, assume that $G$ is a discrete and $A_p(H)$ is Arens regular for every countable subgroup $H$ of $G$. Let $T\in PM_p(G)$. It is enough to show that $\{uT : u\in C_{00}(G), \| u\|_{A_p(G)}\leq 1\}$ is relatively weakly compact in $P M_p(G)$ \cite[Proposition 3]{13}. Now let $u_n$ and $v_m$ be two sequences in the unit ball of $A_p(G)$ with compact support, such that
\begin{equation*}
a =\lim_n \lim_m T (u_nv_m),\qquad	b=\lim_m \lim_n T (u_nv_m)
\end{equation*}
both exist. Assume now that $S = (\cup_{i=1}^{\infty}supp u_i)\cup(\cup_{i=1}^{\infty}=supp v_i)$. Let $H$ be the subgroup of $G$ generated by $S$. Then $H$ is countable. Since $A_p(H)$ is Arens regular. So, $a = b$. It follows that $A_p(G)$ is Arens regular. \\

In \cite{15} Lau and Wong have shown that, if $A_2(G)$ is Arens regular then $G$ is finite, for a locally compact amenable group $G$. Furthermore, Forrest \cite{11} showed that if $G$ is a locally finite discrete group then, $A_p(G)$ is Arens regular if and only if $G$ is finite. In the next theorem we show that the Arens regularity of $A_p(G)$ and finiteness of $G$ can be considered only for finitely generated groups.
\begin{theorem}
Let $G$ be a locally compact group and $A_p(G)$ Arens regular. If the Arens regularity of $A_p(H)$ imply that $H$ is finite, for every finitely generated subgroup $H$ of $G$. Then $G$ is finite.
\end{theorem}
{\bf Proof.}
Suppose $A_p(G)$ is Arens regular, then $G$ is discrete \cite[Theorem 3.2]{10}. Let $H$ be a finitely generated subgroup of $G$. Then by \cite{10} Lemma 3.1, $A_p(H)$ is Arens regular. Hence $H$ is finite. So $G$ is locally finite. It follows from \cite{11} Proposition 3 that $G$ is finite. \\

Let $G$ be a countable discrete group. Then clearly $A_p(G)$ is separable. Hence the $w^*$-topology on the unit ball of $P M_p(G)$, the dual space of $A_p(G)$ is metrizible by \cite{8}, p. 426. So, it is sequentially compact, i.e each sequence in the unit ball of $P M_p(G)$ has a $w^*$-convergent subsequence. Let $Y_p$ denote the unit ball of $P M_p(G)$ and
\begin{equation*}
X_p = span\{\delta_x  : x \in G\},
\end{equation*}
where $\delta_x(f) = f(x)$ for $f\in A_p(G)$, $x \in G$. Suppose $Z_p = X_p\cap Y_p$.
\begin{lemma}
Let $G$ be an amenable discrete group. Then
\begin{equation*}
w^* - clZ_p = Y_p.
\end{equation*}
\end{lemma}
{\bf Proof.}
Since $G$ is amenable, so by \cite[Theorem 5]{13}, $P M_p(G) = CON V_p(G)$ the set of all convolution operators on $L_p(G)$, and the conclusion holds.
\begin{remark}
If $(f_n)$, $(g_m)$ are in the unit ball of $A_p(G)$ and $(\phi_k)$ in $Y_p$. Then clearly there are subsequences $(f^{\prime}_n)$, $(g^{\prime}_n)$ and $(\phi^{\prime}_k)$ of $(f_n)$, $(g_m)$ and $(\phi_k)$ respectively, such that
\begin{equation*}
\lim_n \lim_m \lim_k \phi^{\prime}_k(f^{\prime}_n g^{\prime}_m),\qquad	\lim_m \lim_n \lim_k \phi^{\prime}_k(f^{\prime}_n g^{\prime}_m)
\end{equation*}
both exist.
\end{remark}
We now state the main result of this section.
\begin{theorem}
Let $G$ be an amenable countable discrete group. Then the following are equivalent:
\begin{enumerate}[label=({\roman*})]
\item $A_p(G)$ is Arens regular.
\item For each sequences $(f_n)$, $(g_m)$ in the unit ball of $A_p(G)$, and $(\phi_k)$ in $Z_p$
\begin{equation*}
\lim_n \lim_m \lim_k \phi_k(f_ng_m) = \lim_m \lim_n \lim_k \phi_k(f_ng_m)
\end{equation*}
whenever both exist.
\item For each sequences $(f_n)$, $(g_m)$ in the unit ball of $A_p(G)$ and $\phi_k\in Z_p$. Then
\begin{equation*}
\{\phi_k(f_ng_m) : k > n > m\}− \cap \{\phi_k(f_ng_m) : k > m > n\}^{-}\not=\emptyset.
\end{equation*}
\end{enumerate}
\end{theorem}
{\bf Proof.}
$(i)\Longrightarrow (ii)$: Assume that $(\phi_k)$ is a sequence in $Z_p$. Then there exists a subsequence $({\phi_k}_s )$ of $(\phi_k)$ and $\phi\in Y_p$ such that ${\phi_k}_s \longrightarrow \phi$ in the $w^*$-topology of $P M_p(G)$. Thus we have,
\begin{eqnarray*}
&&a = \lim_n \lim_m \lim_k \phi_k(f_ng_m) = \lim_n \lim_m \lim_s {\phi_k}_s(f_ng_m) = \lim_n \lim_m \phi(f_ng_m), \\
&&b = \lim_m \lim_n \lim_k \phi_k(f_ng_m) =\lim_m \lim_n \lim_s {\phi_k}_s(f_ng_m)=\lim_m \lim_n \phi(f_ng_m).
\end{eqnarray*}
Since $A_p(G)$ is Arens regular, we have $a = b$. \\
$(ii)\Longrightarrow (i)$: Let $\phi\in Y_p$. Then by Lemma 3.3 there is a sequence $(\phi_k)$ in $Z_p$ such that $\phi_k\longrightarrow \phi$ in $w^*$- topology of $P M_p(G)$. Now let $(f_n)$, $(g_m)$ be in the unit ball of $A_p(G)$ such that
\begin{equation*}
a= \lim_n \lim_m \phi(f_ng_m),\qquad b=\lim_m \lim_n \phi(f_ng_m).
\end{equation*}
Then we have
\begin{eqnarray*}
&&a = \lim_n \lim_m \phi(f_ng_m)= \lim_n \lim_m \lim_k \phi_k(f_ng_m), \\
&&b =\lim_m \lim_n \phi(f_ng_m)= \lim_m \lim_n \lim_k \phi_k(f_ng_m).
\end{eqnarray*}
So, $a = b$ and $A_p(G)$ is Arens regular. \\
$(ii)\Longrightarrow (iii)$: Let $(ii)$ holds. Then $(iii)$ follows immediately from Remark 3.4. \\
$(iii)\Longrightarrow (ii)$: It is clear. \\

Let $G$ be a discrete group and $f\in l^1(G)$. Since we can write $f=\sum_{n=1}^{\infty}a_n\chi_{x_n}$. So $f=\sum_{n=1}^{\infty}a_n\chi_{x_n}*\chi^{\vee}_e$ and $f=\sum_{n=1}^{\infty}\|a_n\chi_{x_n}\|_q\|\chi_e\|_p=\sum_{n=1}^{\infty}|a_n|=\|f\|_1$. 		
Whenever, $\chi_x$ denotes the characteristic function $1_{\{x\}}$. Therefore, $f\in A_p(G)$ and $\|f\|_{A_p}\leq \|f\|_1$. Furthermore, it is clear that for each $T\in l^{\infty}(G)$ and $f \in l^1(G)$, $Tf\in C_0(G)$. We now investigate the Arens regularity of $l^1(G)$ under various multiplications and norms. 
\begin{theorem}
\begin{enumerate}[label=({\roman*})]
\item $(l^1(G), . , \|.\|_1)$ is Arens regular. 
\item $(l^1(G), . , \|.\|_1)$ is Arens regular if and only if $G$ is finite. 
\item $(l^1(G), . , \|.\|_{A_p})$ is not normed algebra in general. 
\item $(l^1(G), . , \|.\|_{A_p(G)})$ is Arens regular if and only if $A_p(G)$ is Arens regular.
\end{enumerate}
\end{theorem}
{\bf Proof.} (i). Let $(f_n)$ and $(g_m)$ be two sequences in the unit ball of $l^1(G)$ and
$T \in l^1(G)^*$ such that
\begin{equation*}
\lim_n \lim_m T(f_ng_m),\qquad	\lim_m \lim_n T(f_ng_m)
\end{equation*}
both exist. Since $(l^1(G))_1$ is $w^*$-compact, so $(f_n)$ and $(g_m)$ have subnets $w^*$-converging to some $f$ and $g$ in the unit ball of $l^1(G)$, respectively. Since the original limits will be the same as the limits of the subnets. So we have,
\begin{equation*}
\lim_n \lim_m T(f_ng_m)= \lim_n \lim_m g_m(T f_n) = \lim_n g(T f_n) = \lim_n f_n(T g) = T (fg)
\end{equation*}
and similarly,
\begin{equation*}
\lim_m \lim_n T(f_ng_m)= T(fg).
\end{equation*}
Hence by Theorem 1 of \cite{7} $l^1(G)$ with pointwise multiplication is Arens regular. \\
(ii). See \cite{18}. \\
(iii). Let $G = \mathbb{Z}$. If $V_n=\{−n, −n + 1, ..., 0, 1, ..., n\}$ and $\phi_n(k)=\frac{\chi_{V_n}*\chi^{\vee}_{V_n}(k)}{2n+1}$. Then $\phi_n\in A_p( \mathbb{Z})\cap C_{00}( \mathbb{Z})$, $\|\phi_n\|_{A_p( \mathbb{Z})}=\phi_n(0) = 1$ and
\begin{equation*} 
\phi_n(k)=\left\{ 
\begin{array}{rrr} 
\frac{2n+1-|k|}{2n+1} &  |k|\leq 2n \\
0  &  |k|>2n
\end{array} \right. 
\end{equation*}
But
\begin{equation*}
\phi_2 * \phi^{\vee}_{2}(0) = \phi_2(0)^2 + 2 \sum_{k=1}^{4}\phi_2(k)^2 = 1 +\frac{60}{25}.
\end{equation*}
So $\|\phi_2 * \phi^{\vee}_2\|_{A_p( \mathbb{Z})}>\|\phi_2\|_{A_p( \mathbb{Z})} \|\phi_2\|_{A_p( \mathbb{Z})}$. Therefore $(l^1( \mathbb{Z}),*,\|.\|_{A_p( \mathbb{Z})})$ is not a normed algebra. \\
(iv) Since $l^1(G)$ is normed dense in $A_p(G)$. So, $l^1(G)$ is Arens regular if and only if $A_p(G)$ is Arens regular. \\

Put $B_1$ the unit ball of $(l^1(G), . , \| .\|_{A_p( G)})$.
\begin{example}
Let $G = \mathbb{Z}$ and $\phi_n$ be as above. Then $\phi\in B_1$ and 
$$\phi	= 1 + 2( \sum_{k=1}^{2n}\frac{2n+1-k}{2n+1})= 1 + 2\frac{2n(2n + 1)}{2(2n + 1)}= 1 + 2n.$$	
This shows that $B_1$ is unbounded with $\|.\|_1$-norm.
\end{example}
In the following we study the Arens regularity of $A_p(G)$, whenever $B_1$ is bounded under $\|.\|_1$-norm.
\begin{theorem}
Let $B_1$ be bounded under $\|.\|_1$-norm. Then, 
\begin{enumerate}[label=({\roman*})]
\item	$(l^1(G), .,\|.\|_{A_p(G)})$ is Arens regular.
\item	$A_p(G)$ is Arens regular.
\item	$G$ has no infinite Abelian subgroup.
\item	There is no sequence $\{H_n\}$ of finite subgroups, such that
$$|H_1| < |H_2| < |H_3| < . . .$$
\item There is $M > 0$ such that for each $x\in G, O(x) \leq M$, where, $O(x)$ denote the order of element $x\in G$.
\end{enumerate}
\end{theorem}
{\bf Proof.}
(i) Let $\phi\in (l^1(G), .,\|.\|_{A_p(G)})^*$ and $f\in l^1(G)$. Then
$$|\phi(f)|\leq \|\phi\| \|f\|_{A_p(G)}\leq \|\phi\| \|f\|_1,$$
so $\phi\in l^{\infty}(G)$. Now suppose $(f_n), (g_m)\in B_1$ and
$$\lim_n \lim_m \phi(f_ng_m),\qquad	\lim_m \lim_n \phi(f_ng_m)$$
both exist. Since $B_1$ is bounded under $\|.\|_1$-norm. Then there exists $\alpha> 0$ such that $\| f_n\|\leq \alpha$, $\| g_m\|\leq \alpha$ for all $m, n$. Since $(l^1(G), ., \| .\|_1)$ is Arens regular we have
$$\lim_n \lim_m \phi(f_ng_m) = \lim_m \lim_n \phi(f_ng_m).$$
Hence $(l^1(G), ., \|.\|_{A_p(G)})$ is Arens regular. \\
(ii)	By $3.6 (iv)$. \\
(iii)	Let $H$ be an infinite Abelian subgroup, and $x_1, x_2, . . .$ , be a sequence of
distinct element of $H$. Since $A_p(G)$ is Arens regular, so $G$ is a periodic group. Hence each finitely generated subgroup $H$ is finite. Now let $H_n$ be a subgroup of $H$ generated by $x_1, x_2, x_3, . . . , x_n$ and $\phi_n =\frac{\chi_{H_n}*\chi^{\vee}_{H_n}}{|H_n|}$. Then $\phi_n(e) = 1 =\|\phi\|_{A_p(G)}$ and	
\begin{equation*} 
\phi_n(x)=\left\{ 
\begin{array}{rrr} 
1 &  x\in H_n \\
0  &  x\not\in H_n
\end{array} \right. 
\end{equation*}			
So  $\|\phi_n\|_1= |H_n| \geq n$. This shows that $B_1$ is not bounded under $\|.\|_1$-norm. \\
(iv) If $G$ has a sequence $\{H_n\}$ of finite subgroups, such that
$$|H_1| < |H_2| < |H_3| < . . . .$$
Suppose $\phi_n =\frac{\chi_{H_n}*\chi^{\vee}_{H_n}}{|H_n|}$	, then $\phi_n\in B_1$	and $\|\phi_n\|_1\geq |H_n|$ which is a contradiction. \\	
(v) This follows from (iv). 
\begin{corollary}
Let $G$ be locally finite and $B_1$ be bounded under $\|.\|_1$-norm. Then $G$ is finite.		\end{corollary}		
{\bf Proof.} Let $G$ be an infinite locally finite group. Then $G$ has an infinite Abelian subgroup \cite{14}. It follows from 3.7 (iii) that $G$ must be finite. \\

{\bf Acknowledgment.} We would like to thank the Banach algebra center of Excellence for Mathematics, Univercity of Isfahan.

\bibliographystyle{amsplain}

\authorsaddresses{ Department of Mathematics, University of Isfahan, 81745 Isfahan,Iran 
E-mail address: ghaeidamini@mau.ac.ir} 

\authorsaddresses{Department of Mathematics, University of Isfahan, 81745 Isfahan,Iran 
E-mail address: rejali@sci.ui.ac.ir}

\end{document}